\documentclass[12pt]{article}
 \usepackage{latexsym}
 \usepackage{amssymb}
 \usepackage{graphicx}
 \usepackage{amsmath}

 \newtheorem{Theorem}{Theorem}
\newtheorem{Definition}{Definition}
\newtheorem{Proposition}{Proposition}

\newtheorem{Lemma}{Lemma}
\newtheorem{Corollary}{Corollary}

\newcommand{\A}{{\cal A}}

\newcommand{\uu}{{\bf u}}
\newcommand{\vv}{{\bf v}}
\newcommand{\x}{{\bf x}}
\newcommand{\y}{{\bf y}}
\newcommand{\z}{{\bf z}}
\newcommand{\cc}{{\bf c}}

\newcommand{\0}{{\bf 0}}

\newcommand{\w}{{\bf w}}

\newcommand{\bfalpha}{{\boldsymbol{\alpha}}}
\newcommand{\bfbeta}{{\boldsymbol{\beta}}}
\newcommand{\bfgamma}{{\boldsymbol{\gamma}}}

\newcommand{\qed}{\nobreak \ifvmode \relax \else
      \ifdim\lastskip<1.5em \hskip-\lastskip
      \hskip1.5em plus0em minus0.5em \fi \nobreak
      \vrule height0.75em width0.5em depth0.25em\fi}

\def \ep{\hbox{ }\hfill$\Box$}

\begin{document}
\title{The positive semi-definite cone and sum-of-squares cone of Hankel form}

\author{Zhongming Chen
\thanks{ School of Mathematical
Sciences and LPMC, Nankai University, Tianjin 300071, P.R. China.
Email: czm183015@mail.nankai.edu.cn.
This author's work was partially done
when he was visiting The Hong Kong Polytechnic University.}, \quad
Liqun Qi
\thanks{ Department of Applied
Mathematics, The Hong Kong Polytechnic University, Hung Hom,
Kowloon, Hong Kong. Email: maqilq@polyu.edu.hk. This author's work was supported by the Hong
Kong Research Grant Council (Grant No. PolyU 502111, 501212, 501913 and 15302114).}
}

\date{\today}
\maketitle

\begin{abstract}
\noindent
In this paper, the geometry properties of Hankel form are studied,
including their positive semi-definite (PSD) cone and sum-of-squares (SOS) cone.
We denote them by $HPSD(m,n)$ and $HSOS(m,n)$, respectively.
We show that both $HPSD(m,n)$ and $HSOS(m,n)$ are closed convex cones.
The dual cone of $HPSD(m,n)$ is the convex hull of all $m$-times convolutions of real vectors.
Besides, we derive the dual cone of SOS tensors.
By reformulation, it follows that the dual cone of $HSOS(m,n)$ can also be written explicitly.
These results may lead further research on the Hilbert-Hankel problem.
\vspace{3mm}

\noindent {\bf Key words:}\hspace{2mm}
Hankel form, positive semi-definite cone, SOS cone, dual cone
\vspace{2mm}

\noindent {\bf AMS subject classifications (2010):}\hspace{2mm}
15A18; 15A69

\end{abstract}

\newpage
\section{Introduction}
\hspace{14pt}
Hankel structures are frequently encountered in applications such as signal processing \cite{o2001}.
Besides Hankel matrices, tensors with different Hankel structures also find applications in
high-order singular value decompositions (HOSVD) \cite{bb2008}, exponential data fitting
\cite{dqw2015,pdv2005,pdv2009} and signal separation \cite{d2011}.
In this paper, we study the geometry propositions of Hankel form, including
their positive semi-definite (PSD) cone and sum-of-squares (SOS) cone.

Let $\A= (a_{i_1 \cdots i_m}) \in T_{m,n}$.
If there is a vector $\vv = (v_0, v_1, \cdots, v_{(n-1)m})^\top \in \mathbb{R}^{(n-1)m+1}$
such that for $i_1, \cdots, i_m  \in [n]$, we have
$$  a_{i_1 \cdots i_m} \equiv v_{i_1 +i_2+ \cdots + i_m - m},   $$
then we say that $\A$ is an $m$th order {\bf Hankel tensor} and the vector
$\vv$ is called the { generating vector} of $\A$ \cite{clq20151,cqw20151,cqw2015,dqw2015,lqw20141,lqx2014,q2015,x2015}. Clearly, a Hankel tensor is symmetric.
Let $\x=(x_0,x_1,\cdots, x_p)^\top \in \mathbb{R}^{p+1}$ and $\y=(y_0,y_1,\cdots, y_q)^\top \in \mathbb{R}^{q+1}$.
The {\bf convolution} of $\x$ and $\y$ is defined as $\z = \x * \y =(z_0,z_1, \cdots, z_{p+q}) \in \mathbb{R}^{p+q+1}$, with
$$  z_i = \sum_{ i_1 + i_2 =i \atop 0 \leq i_1 \leq p, 0 \leq i_2 \leq q  } x_{i_1} y_{i_2} $$
for $i =0, 1, \cdots, p+q$ \cite{dqw20152}. If we regard the entries of $\x$ and $\y$ as the coefficients of two polynomials, then the entries
of $\z$ is the coefficients of the product polynomial of these two polynomials. For any $\x \in \mathbb{R}^n$, we denote
$$ \x^{*m} = \overbrace{\x * \cdots  * \x}^m . $$
Let $\vv \in \mathbb{R}^{(n-1)m+1}$ and $\x \in \mathbb{R}^n$. Then
$ H = \vv \bullet \x^{*m}  $
is called a {\bf Hankel form} of order $m$ and dimension $n$, where '$\bullet$' denotes the standard inner product.
In fact, let $\A \in T_{m,n}$ be the Hankel tensor generated by $\vv$, it follows that
$ H = \vv \bullet \x^{*m} = \A \x^m $.

Let $m$ be even. We say that a Hankel form $H$ is PSD if for any $\x \in \mathbb{R}^n$, $H \geq 0$.
We say that $H$ is SOS if $H$ can be written as a sum of squares of homogeneous polynomials
of $\x$ with degree $k$. Clearly, if $H$ is SOS, then $H$ is PSD. The main question is:
if $H$ is PSD, is it SOS? If the answer to this question is yes, then the problem for determining an even order Hankel tensor
is positive semi-definite or not is solvable in polynomial-time \cite{lqx2014}.
In a certain sense, it is the 17th Hilbert problem with the Hankel constraint.
And the problem raised by the above question is called the Hilbert-Hankel problem \cite{cqw2015}.
Recently, some work showed that there are no PSD non-SOS Hankel tensors under certain conditions \cite{clq20152,cqw2015,lqw20141,lqw20142,lqx2014,q2015}.
However, until now, it still remains an open problem.

In this paper, we study the PSD cone and SOS cone of Hankel form respectively.
This may lead further research on the Hilbert-Hankel problem.
The rest of this chapter is organized as follows.
In Section 2, we consider the PSD cone of Hankel form.
And we show that its dual cone is the convex hull of all $m$-times convolutions of real vectors.
The Hankel spectrahedra is introduced in Section 3. By reformulation,
we give a new way to get the dual cone of the PSD cone and SOS cone of Hankel form.
In Section 4, we study the SOS cone of Hankel form, as well as its dual cone.

\section{The PSD cone of Hankel form}
\hspace{14pt}
Suppose that $m=2k$, $\vv=(v_0,\cdots,v_{(n-1)m})^\top$ and $\x=(x_0,\cdots,x_{n-1})^\top$. Define
$$ HPSD(m,n)=\{\vv \in \Re^{(n-1)m+1} : H=\vv \bullet \x^{*2k} \text{ is PSD}\}   $$
and
$$ HSOS(m,n)=\{\vv \in \Re^{(n-1)m+1} : H=\vv \bullet \x^{*2k} \text{ is SOS}\}.  $$
Obviously, $HSOS(m,n) \subseteq HPSD(m,n)$.
\begin{Proposition}
Let $m$ be even. Then $HPSD(m,n)$ and $HSOS(m,n)$ are closed convex cone.
\end{Proposition}
{\bf Proof.}
Clearly, $HPSD(m,n)$ is a convex cone. Now we prove its closeness.
Let $\{ \vv_k\} \subseteq HPSD(m,n)$ with $\vv_k \rightarrow \vv$.
By definition, we have $\vv_k \bullet \x^{*m} \geq 0$ for any $\x \in \Re^n$.
It follows that $\vv \bullet \x^{*m} = \lim_{k \rightarrow \infty} \vv_k \bullet \x^{*m} \geq 0$
for any $\x \in \Re^n$. Hence, $\vv \in HPSD(m,n)$.
The convexity of the $HSOS(m,n)$ cone also directly follows from the definition.
Let $\{ \vv_k\} \subseteq HSOS(m,n)$ with $\vv_k \rightarrow \vv$. Since $HSOS(m,n) \subseteq HPSD(m,n)$,
$\vv \in HPSD(m,n)$. For $\x \in \Re^n$, $\vv_k \bullet \x^{*m}$ is an SOS polynomial and
$\vv_k \bullet \x^{*m} \rightarrow \vv \bullet \x^{*m}$. Note from \cite{Lau} that the set of all SOS
polynomials on $\Re^n$ with degree at most $m$ is a closed cone.
So, $\vv \bullet \x^{*m}$ is also an SOS polynomial. Therefore, $\vv \in HSOS(m,n)$.
\ep

\medskip
Recall that for a given set $S$ in the Euclidean space $\Re^n$, its dual cone $S^*$ is defined as
$$
S^* = \{ \y \in \Re^n : \y \bullet \x \geq 0 \text { for all } \x \in S  \}.
$$
Note that $S^*$ is always a convex cone, even if $S$ is neither convex nor a cone.
When $S$ is a closed convex cone, we have $S^{**} =S$. To establish the dual cone of $HPSD(m,n)$,
we introduce a set which is the convex hull of all $m$-times convolutions of real vectors.

\begin{Definition}\label{d1}
Let $m$ be even and $n \in \mathbb{N}$. We denote by $U(m,n)$ the convex hull of all $m$-times convolutions of real vectors in $\Re^n$, i.e.,
$$ U(m,n):={\rm conv}\{ \x^{*m} : \x \in \Re^n \} .$$
\end{Definition}

In fact, $U(m,n)$ is a closed convex cone. Before that, a useful lemma is needed.
\begin{Lemma}
Let $m$ be even and $n \in \mathbb{N}$. For any $\x \in \Re^n$ and $\y \in \Re^n$, there
exists a constant $c >0$ such that
$$ c \| \x^{*m} + \y^{*m} \| \geq  \max\{ \| \x \|^m, \| \y \|^m \} . $$
\end{Lemma}
{\bf Proof.}
Let $\z = \x^{*m} + \y^{*m}$. For any $t \in \Re$, denote
$$ {\bf t} =(1, t, \cdots, t^{n-1})^\top \quad  \text{ and }  \quad {\bf \tilde{t}} = (1, t, \cdots, t^{(n-1)m})^\top .$$
By definition, we have
$\z \bullet {\bf \tilde{t}} = (\x \bullet {\bf t})^m + (\y \bullet {\bf t})^m \geq \max\{ (\x \bullet {\bf t})^m , (\y \bullet {\bf t})^m\}$.
According to Cauchy-Schwarz inequality, one can obtain that $\| \z \| \| {\bf \tilde{t}} \| \geq \max\{ (\x \bullet {\bf t})^m , (\y \bullet {\bf t})^m\}$.

Now for $k=0, \cdots, n-1$, we choose $t=t_k \in \Re$ such that $t_i \neq t_j$ when $i\neq j$. Denote
$$ {\bf t}_k =(1, t_k, \cdots, t_k^{n-1})^\top  \quad  \text{ and }  \quad {\bf \tilde{t}}_k = (1, t_k, \cdots, t_k^{(n-1)m})^\top .$$
Let $c_k = \| {\bf \tilde{t}}_k \| >0$. Define $\w =(w_0, \cdots, w_{n-1})^\top \in \Re^n $ with $w_k = \x \bullet {\bf t}_k$.
Then we have $\w = T \x$, where $T$ is the transpose of a Vandermonde matrix with the $k$th column ${\bf t}_k$, $k=0,\cdots, n-1$.
It follows that
$$ \| \z \| \sum_{k=1}^{n-1} c_k  \geq \sum_{k=0}^{n-1} w_k^m  = \| \w \|_m^m.  $$
On the other hand, there exists a constant $C >0$ such that $\| \w \| \leq C \| \w \|_m$.
Since $T$ is nonsingular, we have $\| \x \| \leq \| T^{-1} \|  \| \w \|$, where $\| T^{-1} \| >0$ is the largest singular of $T^{-1}$.
So $$ \| \x \|^m \leq \| T^{-1} \|^m \sum_{k=0}^{n-1} c_k \| \z \| .$$
Let $c= \| T^{-1} \|^m \sum_{k=0}^{n-1} c_k >0$. Then $\| \x \|^m \leq c \| \z \| $.
By a similar way, we have $ \| \y \|^m \leq c \| \z \| $. So the proof is completed.
\ep

\medskip
From the proof above, we have the following corollary.
\begin{Corollary}\label{c1}
Let $m$ be even and $n \in \mathbb{N}$. Let $s \geq 1$ be a given integer. Then for any $\x^j \in \Re^n$, $j=1,\cdots,s$, there exists a constant $c>0$
such that $$c \| \sum_{j=1}^s (\x^j)^{*m} \| \geq  \max_{1 \leq j \leq s} \| \x^j \|^m .$$
\end{Corollary}

We are now ready to prove that $U(m,n)$ is a closed convex cone.
\begin{Lemma}\label{l2}
Let $m$ be even and $n \in \mathbb{N}$. Then $U(m,n)$ is a closed convex cone with dimension at most $(n-1)m+1$.
\end{Lemma}
{\bf Proof.} Since $U(m,n) \subseteq \Re^{(n-1)m+1}$, the dimension of $U(m,n)$ is at most $(n-1)m+1$.
The convexity of the $U(m,n)$ cone also directly follows from the definition. To see the closeness,
let $\{ \y_k \}_{k=1}^{\infty} \subseteq U(m,n)$ such that $ \y_k \rightarrow \y$. Clearly, $\y_k$ is bounded.
For any $\y_k$, by the Carath\'{e}odory theorem, there exist $\x_k^j$, $j=0,\cdots,(n-1)m$ such that
$$ \y_k = \sum_{j=0}^{(n-1)m} (\x_k^j)^{*m}  .$$
By Corollary \ref{c1}, we can see that for any $j=0,\cdots, (n-1)m$, the sequence $\{ \x_k^j\}_{k=1}^{\infty}$
is also bounded. By passing to the subsequences, we can assume that $\x_k^j \rightarrow \x^j $, $j=0,\cdots,(n-1)m$.
It follows that $$\y = \lim_{k \rightarrow \infty} \y_k = \sum_{j=0}^{(n-1)m} (\x^j)^{*m} \in U(m,n). $$
Thus, the proof is completed.
\ep

\begin{Theorem}
Let $m$ be even and $n \in \mathbb{N}$. Then
$$ U(m,n)^* = HPSD(m,n)  \quad \text{ and } \quad HPSD(m,n)^* = U(m,n) .$$
\end{Theorem}
{\bf Proof.}
First, we show that $U(m,n)^* = HPSD(m,n)$.
Suppose $\vv \in U(m,n)^*$. Then, by definition, we have $\vv \bullet \y \geq 0 $ for all $\y \in U(m,n)$.
In particular, $\vv \bullet \x^{*m} \geq 0$ for all $\x \in \Re^n$. It follows that $\vv \in HPSD(m,n)$.
Hence, $U(m,n)^* \subseteq HPSD(m,n)$. On the other hand, the dimension of $U(m,n)$ is at most $(n-1)m+1$.
For any $\y \in U(m,n)$, by the Carath\'{e}odory theorem, there exist $\x^j \in \Re^n$, $j=0,\cdots (n-1)m$
such that $$\y = \sum_{j=0}^{(n-1)m} (\x^j)^{*m}.$$
Suppose $\vv \in HPSD(m,n)$. Then $\vv \bullet \y \geq 0$ since $\vv \bullet (\x^j)^{*m} \geq 0 $ for any $j=0,\cdots (n-1)m$.
So $HPSD(m,n) \subseteq U(m,n)^*$. Therefore, we have $U(m,n)^* = HPSD(m,n)$.

For the second part, the equality $ HPSD(m,n)^* = U(m,n)^{**} = U(m,n)$ holds since $U(m,n)$ is a closed convex cone by Lemma \ref{l2}.
\ep

\medskip
According to the fact $HSOS(m,n) \subseteq HPSD(m,n)$, we have the following corollary.
\begin{Corollary}
Let $m$ be even and $n \in \mathbb{N}$. Then
$$ HSOS(m,n) \subseteq U(m,n)^* \quad \text{ and } \quad U(m,n) \subseteq HSOS(m,n)^*    . $$
\end{Corollary}

\section{Hankel Spectrahedra}
\hspace{14pt}
As we know, Semidefinite programming is a very important and valuable tool in polynomial optimization.
And it can be solved in polynomial time.
The feasible region of a Semidefinite Program forms a closed convex set which has the form
$$
S= \{ \x \in \Re^m : Q_0 + \sum_{i=1}^m x_i Q_i \succeq 0  \},
$$
where $Q_i \in S_{2,n}$, $i=0,1,\cdots,m$.
The feasible region of a Semidefinite Program is called a {\bf Spectrahedra} \cite{RG}.
When $Q_0 = 0$, $S$ is a closed convex cone.

Given a set of Vandermonde vectors
\begin{equation}\label{e0}
\uu_k = (1, u_k,\cdots, u_k^{n-1})^\top \in \Re^n, \quad k=0,\cdots,(n-1)m,
\end{equation}
where $u_i \neq u_j$ when $i \neq j$.
It has been shown \cite{q2015} that a symmetric tensor $\A \in S_{m,n}$ is a Hankel tensor
if and only if it has a Vandermonde decomposition, i.e., there exists a vector ${\bfalpha} =(\alpha_0, \cdots, \alpha_{(n-1)m})^\top \in \Re^{(n-1)m+1}$
such that
$$
\A = \sum_{k=0}^{(n-1)m} \alpha_k \uu_k^{\otimes m},
$$
where $\uu^{\otimes m} \in S_{m,n}$ is the rank-one tensor with entries $u_{i_1} \cdots u_{i_m}$.
In fact, this representation of a Hankel tensor is unique when the set of Vandermonde vectors is given.

Let $m$ be even. A tensor $\A \in S_{m,n}$ is called positive semi-definite if $\A \x^m \geq 0$ for all $\x \in \Re^n$.
Let $PSD_{m,n}$ be the set of all positive semi-definite tensors. Denote
$$
G = \{ \bfalpha \in \Re^{(n-1)m+1} : \sum_{k=0}^{(n-1)m} \alpha_k \uu_k^{\otimes m} \in PSD_{m,n} \}.
$$
This can be seen as the generalized spectrahedra for the tensor case.
Clearly, $\Re^{(n-1)m+1}_+ \subseteq G$. It is easy to see that $G$ is also a closed convex cone.
And there is a one-to-one mapping between $G$ and $HPSD(m,n)$.
In particular, let $U \in \Re^{[(n-1)m+1] \times [(n-1)m+1]}$ be the Vandermonde matrix with the $k$th
column $(1, u_k,\cdots, u_k^{(n-1)m})^\top$, $k = 0,1.\cdots, (n-1)m$.
Since $u_i \neq u_j$ when $i \neq j$, the Vandermonde matrix is nonsingular. And the linear mapping
$f: G \rightarrow HPSD(m,n)$ is a bijection defined by
$f(\bfalpha) = U \bfalpha $.

By this reformulation, we give a new way to get the dual cone of $HPSD(m,n)$.
Denote by $U_{m,n}$ the
convex hull of all $m$th-order n-dimensional symmetric rank-one tensors, i.e.,
$$
U_{m,n} = \text{conv} \{ \x^{\otimes m} : \x \in \Re^n  \}.
$$
It has been shown \cite{HLQ} that $PSD_{m,n}$ and $U_{m,n}$ are dual cones, i.e.,
$$
PSD_{m,n}^* =U_{m,n}   \quad \text{ and }  \quad  U_{m,n}^* = PSD_{m,n}.
$$
Due to this result, we have the following conclusion.
\begin{Theorem}\label{t2}
Let $G$ be defined as above. The dual cone of $G$ is the closure of the set
$$ H: = \left\{ (\A \uu_0^m, \A \uu_1^m, \cdots, \A \uu_{(n-1)m}^m )^\top  \in \Re^{(n-1)m+1} :  \A \in U_{m,n}  \right\}.  $$
\end{Theorem}
{\bf Proof.} Clearly, $H$ is convex cone. Suppose $\bfalpha \in H^*$. By definition, we have that
$$ \sum_{k=0}^{(n-1)m} \alpha_k \A \uu_k^m = \left\langle \A , \sum_{k=0}^{(n-1)m} \alpha_k \uu_k^{\otimes m}  \right\rangle  \geq 0, \quad \forall \A \in U_{m,n} .$$
It means that $\bfalpha \in H^*$ if and only if $\sum_{k=0}^{(n-1)m} \alpha_k \uu_k^{\otimes m} \in U_{m,n}^* = PSD_{m,n}$.
So we have $H^*=G$.
Note that $H$ is not closed in general. It then follows from the double polar theorem in convex analysis
that $G^* =H^{**}= \text{cl}H$.
\ep

\begin{Corollary}
Let $m$ be even and $n \in \mathbb{N}$. Let $U \in \Re^{[(n-1)m+1] \times [(n-1)m+1]}$ be the Vandermonde matrix generated by the set $\{ u_k\}_{k=0}^{(n-1)m}$
with $u_i \neq u_j$ when $i \neq j$. Then
$$ \{ U^\top \y : \y \in U(m,n) \} =  \text{\rm cl}H ,$$
where $U(m,n)$ and $H$ are defined in Definition \ref{d1} and Theorem \ref{t2}, respectively.
\end{Corollary}

\section{The SOS cone of Hankel form}
\hspace{14pt}
Let $m$ be even. A tensor $\A \in S_{m,n}$ is called a SOS tensor if $\A \x^m \geq 0$ is a SOS polynomial of $\x \in \Re^n $ \cite{lqy2015}.
Let $SOS_{m,n}$ be the set of all SOS tensors. Denote
$$
SG = \{ \bfalpha \in \Re^{(n-1)m+1} : \sum_{k=0}^{(n-1)m} \alpha_k \uu_k^{\otimes m} \in SOS_{m,n} \},
$$
where $\uu_k \in \Re^n$ are Vandermonde vectors given in (\ref{e0}).
Clearly, $\Re^{(n-1)m+1}_+ \subseteq SG$.
From the analysis above, the linear mapping $\tilde{f}: SG \rightarrow HSOS(m,n)$ defined by
$\tilde{f}(\bfalpha) = U \bfalpha $ is a one-to-one mapping between $SG$ and $HSOS(m,n)$,
where $U \in \Re^{[(n-1)m+1] \times [(n-1)m+1]}$ is the Vandermonde matrix with the $k$th
column $(1, u_k,\cdots, u_k^{(n-1)m})^\top$, $k = 0,1.\cdots, (n-1)m$.
Since $u_i \neq u_j$ when $i \neq j$, the Vandermonde matrix is nonsingular.

In fact, if we know the dual cone of $SOS_{m,n}$, the dual cone of $HSOS(m,n)$ can also be derived similarly.
First, we consider the SOS cone of symmetric tensors with $(m,n)=(6,3)$, i.e., $SOS_{6,3}$.
For any $\A \in S_{6,3}$, there are 28 independent elements.
We index these elements by the $3$-tuples of degree $6$, i.e.
$$
\A = (a_{\bfalpha}), \quad  \bfalpha \in \left\{ (\alpha_0, \alpha_1, \alpha_2) : \sum_{i=0}^2 \alpha_i = 6, \alpha_i \geq 0, \forall i=0,1,2  \right\}.
$$
For simplicity, we write $\A =(a_{\bfalpha})_{|\bfalpha|=6} \in S_{6,3}$ to emphasis on these 28 independent elements.
For any $\w =(x, y,z)^\top \in \Re^3$, $\A \w^6$ is a homogeneous polynomial with degree 6, i.e.,
$$ \A \w^6 = \sum_{ | \bfalpha| =6} c_{\bfalpha} a_{\bfalpha} \w^{\bfalpha} ,$$
where $\w^{\bfalpha}$ denotes the monomial $x^{\alpha_0} y^{\alpha_1} z^{\alpha_2}$, and $c_{\bfalpha}$ is the number of the elements $a_{\bfalpha}$
in $\A$. For example, $c_{600}=1$, $a_{510} = 6$ and $c_{321} = 60$. By simple computation, there are exactly 28 monomials of $\w$ of degree 6.

Let $[\w]_3 := (x^3, x^2y, x^2z, xy^2, xyz, y^3, xz^2, y^2z, yz^2,z^3 )^\top \in \Re^{10}$ be the vector of all the monomials of $\w$ of degree 3.
Then $\A \w^6 $ is a sum of square with SOS rank $r$ if and only if there exist $\cc_i \in \Re^{10}\setminus \{\0 \}$, $i=1,\cdots, r$, such that
$$ \A \w^6 = \sum_{i=1}^r \left(  \cc_i^\top [\w]_3 \right)^2 = [\w]_3^\top \left( \sum_{i=1}^r \cc_i \cc_i^\top \right) [\w]_3 , $$
i.e., there exists a positive semi-definite matrix $Q \in \Re^{10 \times 10}_+ $ with rank $r$ such that
$$ \A \w^6 =  [\w]_3^\top Q [\w]_3 .$$
Moreover, by indexing the matrix $Q$ by the 10 monomials of $\w$ of degree 3
(or, more precisely, the associated exponent tuples), we obtain the following conditions:
$$
c_{\bfalpha} a_{\bfalpha} = \sum_{\bfbeta + \bfgamma = \bfalpha} Q_{\bfbeta  \bfgamma}, \quad  Q \succeq 0 .
$$
This is a system of 28 linear equations, one for each coefficient of $\A \w^6$.
We denote these linear equations by $ \langle A_{\bfalpha} , Q \rangle = c_{\bfalpha} a_{\bfalpha}$ for any $\bfalpha$ satisfying $|\bfalpha | =6. $
In particular, the matrix $Q$ is expressed as follows:
\begin{align}
\begin{bmatrix}\begin{smallmatrix}
q_{(300,300)} & q_{(300,210)} & q_{(300,201)} & q_{(300,120)} & q_{(300,111)} & q_{(300,030)} & q_{(300,102)} & q_{(300,021)} & q_{(300,012)} & q_{(300,003)}  \\
q_{(210,300)} & q_{(210,210)} & q_{(210,201)} & q_{(210,120)} & q_{(210,111)} & q_{(210,030)} & q_{(210,102)} & q_{(210,021)} & q_{(210,012)} & q_{(210,003)}  \\
q_{(201,300)} & q_{(201,210)} & q_{(201,201)} & q_{(201,120)} & q_{(201,111)} & q_{(201,030)} & q_{(201,102)} & q_{(201,021)} & q_{(201,012)} & q_{(201,003)}  \\
q_{(120,300)} & q_{(120,210)} & q_{(120,201)} & q_{(120,120)} & q_{(120,111)} & q_{(120,030)} & q_{(120,102)} & q_{(120,021)} & q_{(120,012)} & q_{(120,003)}  \\
q_{(111,300)} & q_{(111,210)} & q_{(111,201)} & q_{(111,120)} & q_{(111,111)} & q_{(111,030)} & q_{(111,102)} & q_{(111,021)} & q_{(111,012)} & q_{(111,003)}  \\
q_{(030,300)} & q_{(030,210)} & q_{(030,201)} & q_{(030,120)} & q_{(030,111)} & q_{(030,030)} & q_{(030,102)} & q_{(030,021)} & q_{(030,012)} & q_{(030,003)}  \\
q_{(102,300)} & q_{(102,210)} & q_{(102,201)} & q_{(102,120)} & q_{(102,111)} & q_{(102,030)} & q_{(102,102)} & q_{(102,021)} & q_{(102,012)} & q_{(102,003)}  \\
q_{(021,300)} & q_{(021,210)} & q_{(021,201)} & q_{(021,120)} & q_{(021,111)} & q_{(021,030)} & q_{(021,102)} & q_{(021,021)} & q_{(021,012)} & q_{(021,003)}  \\
q_{(012,300)} & q_{(012,210)} & q_{(012,201)} & q_{(012,120)} & q_{(012,111)} & q_{(012,030)} & q_{(012,102)} & q_{(012,021)} & q_{(012,012)} & q_{(012,003)}  \\
q_{(003,300)} & q_{(003,210)} & q_{(003,201)} & q_{(003,120)} & q_{(003,111)} & q_{(003,030)} & q_{(003,102)} & q_{(003,021)} & q_{(003,012)} & q_{(003,003)}
\end{smallmatrix}
\end{bmatrix}\nonumber
\end{align}
Let $E_{ij}$ be the matrix in $\Re^{10 \times 10}$ with $(i,j)$-th entry 1 and 0 otherwise.
It follows that
$A_{600} = E_{11},  A_{510} = E_{12}+E_{21},  A_{501} =E_{13}+E_{31},  A_{420} = E_{14} + E_{22}+E_{41}$ and so on.
We can see that $\langle A_{\bfalpha}, A_{\bfbeta} \rangle =0$ if $\bfalpha \neq \bfbeta$.
By this reformulation, we have
\begin{equation}\label{e1}
SOS_{6,3} = \left\{ (a_{\bfalpha})_{|\bfalpha|=6} \in S_{6,3} : \langle A_{\bfalpha}, Q \rangle = c_{\bfalpha} a_{\bfalpha} \text{ for all } |\bfalpha| =6, \ \   Q \in \Re^{10 \times 10}_+ \right\} .
\end{equation}

In fact, we can generalize this result to $SOS_{m,n}$ where $m$ is even. Let $m=2k$.
For any $\x =(x_0, \cdots, x_{n-1})^\top \in \Re^n$, there are $\binom{n+m-1}{m}$ monomials of $\x$ of degree $m$,
and there are $\binom{n+k-1}{k}$ monomials of $\x$ of degree $k$.
It means that for $SOS_{m,n}$, there are $\binom{n+m-1}{m}$ linear equations on the right hand side of (\ref{e1}),
and $Q \in \Re^{\binom{n+k-1}{k} \times \binom{n+k-1}{k}}_+$. Moreover, we can write the constant $c_{\bfalpha}$ precisely.
Denote $S = \left\{ (\alpha_0, \cdots, \alpha_{n-1}) : \sum_{i=0}^{n-1} \alpha_i =m , \alpha_i \geq 0, \forall i=0,\cdots, n-1 \right\}$.
For any tuple of exponents $\bfalpha \in S $,
we have
\begin{equation}\label{e2}
c_{\bfalpha} = \binom{m}{\alpha_0} \binom{m-\alpha_0}{\alpha_1} \cdots \binom{\alpha_{n-1}}{\alpha_{n-1}} .
\end{equation}
For instance, $c_{m0\cdots 0} = 1$, $c_{(m-1) 1 \cdots 0} = m $ and $c_{00\cdots m} = 1$.
\begin{Theorem}\label{t3}
Let $m = 2k$ be even and $n \in \mathbb{N}$. Then
$$
SOS_{m,n} = \left\{ (a_{\bfalpha})_{|\bfalpha|=m} \in S_{m,n} : \langle A_{\bfalpha}, Q \rangle = c_{\bfalpha} a_{\bfalpha} \text{ \rm for all } |\bfalpha| =m,  \ \  Q  \in \Re^{\binom{n+k-1}{k} \times \binom{n+k-1}{k}}_+ \right\},
$$
where $c_{\bfalpha}$ are given by (\ref{e2}) and the $\binom{n+m-1}{m}$ matrices $\{ A_{\bfalpha} \}_{|\bfalpha|=m}$ are only depend on $(m,n)$, satisfying
$\langle A_{\bfalpha}, A_{\bfbeta} \rangle =0$ if $\bfalpha \neq \bfbeta$.
\end{Theorem}

Based on this result, the dual cone of $SOS_{m,n}$ can be derived immediately.

\begin{Theorem}\label{t4}
Let $m = 2k$ be even and $n \in \mathbb{N}$. Then,
$$  SOS_{m,n}^* = \left\{ (b_{\bfalpha})_{|\bfalpha|=m} \in S_{m,n} : \sum_{|\bfalpha|=m} b_{\bfalpha} A_{\bfalpha} \succeq 0  \right\} , $$
where $\{ A_{\bfalpha} \}_{|\bfalpha|=m} $ are symmetric $\binom{n+k-1}{k} \times \binom{n+k-1}{k}$ matrices given in Theorem \ref{t3}.
\end{Theorem}
{\bf Proof.}
Let $M = \left\{ (b_{\bfalpha})_{|\bfalpha|=m} \in S_{m,n} : \sum_{|\bfalpha|=m} b_{\bfalpha} A_{\bfalpha} \succeq 0  \right\}$.
Let $(b_{\bfalpha})_{|\bfalpha| =m} \in S_{m,n}$. It follows from Theorem \ref{t3} that
$(b_{\bfalpha}) \in SOS_{m,n}^*$ if and only if
$ \langle (b_{\bfalpha}) , (a_{\bfalpha}) \rangle = \sum_{|\bfalpha| =m } c_{\bfalpha} b_{\bfalpha} a_{\bfalpha} \geq 0$
for any $(a_{\bfalpha}) \in SOS_{m,n}$, i.e.,
$$
 \sum_{|\bfalpha| =m } b_{\bfalpha} \left\langle A_{\bfalpha}, Q \right\rangle = \left\langle \sum_{|\bfalpha| =m } b_{\bfalpha} A_{\bfalpha}, Q \right\rangle \geq 0, \quad \forall Q \in \Re^{\binom{n+k-1}{k} \times \binom{n+k-1}{k}}_+ .
$$
Since all the positive semi-definite matrices are self-dual, this is equivalent to say that $ (b_{\bfalpha})_{|\bfalpha| =m} \in M$.
The proof is completed.
\ep

\medskip
It can be seen that the dual cone of $SOS_{m,n}$ is nothing but a spectrahedra.

In the following, we establish the dual cone of $HSOS(m,n)$.
Let $SG$ be the set defined in the beginning of this section.
Similar with Theorem \ref{t2}, we have
\begin{equation}\label{e3}
SG^* = \text{cl}\left\{  (\A \uu_0^m, \A \uu_1^m, \cdots, \A \uu_{(n-1)m}^m )^\top  \in \Re^{(n-1)m+1} :  \A \in SOS_{m,n}^*    \right\} ,
\end{equation}
where $\uu_k \in \Re^n$ are Vandermonde vectors given in (\ref{e0}).
Since there is a one-to-one mapping between $SG$ and $HSOS(m,n)$, we have the following conclusion.
\begin{Theorem}
Let $m$ be even and $n \in \mathbb{N}$. Let $U \in \Re^{[(n-1)m+1] \times [(n-1)m+1]}$ be the Vandermonde matrix generated by the set $\{ u_k\}_{k=0}^{(n-1)m}$
with $u_i \neq u_j$ when $i \neq j$. Then
$$ HSOS(m,n)^* = \left\{ (U^\top)^{-1} \bfalpha : \bfalpha \in SG^* \right\},$$
where $SG^*$ is given by (\ref{e3}).
\end{Theorem}


\medskip
\begin{thebibliography}{abc99xyz}


\bibitem{bb2008} Badeau R, Boyer R. Fast multilinear singular value decomposition for structured tensors.
    SIAM J Matrix Anal Appl, 2008, 30: 1008-1021

\bibitem{clq20151} Chen H, Li G, Qi L. Further results on Cauchy tensors and Hankel tensors. http://arxiv.org/abs/1501.06726, 2015

\bibitem{clq20152} Chen H, Li G, Qi L. Sum-of-squares tensors and their sum-of-squares rank. http://arxiv.org/abs/1504.03414, 2015

\bibitem{cqw20151} Chen Y, Qi L, Wang Q. Computing eigenvalues of large scale Hankel tensors. http://arxiv.org/abs/1504.07413, 2015

\bibitem{cqw2015} Chen Y, Qi L, Wang Q. Positive Semi-Definiteness and Sum-of-Squares Property of Fourth Order Four Dimensional Hankel Tensors.
    http://arxiv.org/abs/1502.04566, 2015

\bibitem{d2011} De Lathauwer L. Blind separation of exponential polynomials and the decomposition of a tensor in rank-$(L_r,L_r,1)$
terms. SIAM J Matrix Anal Appl, 2011, 32: 1451-1474

\bibitem{dqw2015} Ding W, Qi L, Wei Y. Fast Hankel tensor-vector product and its application to exponential data fitting,
to appear in: Numerical Linear Algebra with Applications

\bibitem{dqw20152} Ding W, Qi L, Wei Y. Positive Semi-Definite Hankel Tensors. http://arxiv.org/abs/1505.02528, 2015

\bibitem{HLQ}  Hu S, Li G, Qi L. A tensor analogy of Yuan's alternative theorem and polynomial optimization with sign structure,
to appear in: Journal of Optimization Theory and Applications

\bibitem{Lau} Laurent M. Sums of squares, moment matrices and optimization over polynomials.
Emerging Applications of Algebraic Geometry, Vol. 149 of IMA Volumes
in Mathematics and its Applications, M. Putinar and S. Sullivant eds., Springer,
2009, 157-270.

\bibitem{lqw20141} Li G, Qi L, Wang Q. Are there sixth order three dimensional PNS Hankel tensors?  http://arxiv.org/abs/1411.2368, 2014

\bibitem{lqw20142} Li G, Qi L, Wang Q. Positive semi-definiteness of generalized anti-circulant tensors, to appear in:
Communications in Mathematical Sciences

\bibitem{lqx2014}  Li G, Qi L, Xu Y. SOS Hankel Tensors: Theory and Application. http://arxiv.org/abs/1410.6989, 2014

\bibitem{lqy2015} Luo Z, Qi L, Ye Y. Linear operators and positive semidefiniteness of symmetric tensor spaces.
Sci China Math, 2015, 58: 197-212

\bibitem{o2001} Olshevsky V. Structured Matrices in Mathematics, Computer Science, and Engineering. II, Contemporary Mathematics, vol. 281.
American Mathematical Society: Providence, RI, 2001

\bibitem{pdv2005} Papy J M, De Lathauwer L, Van Huffel S. Exponential data fitting using multilinear algebra: the single-channel and
multi-channel case. Numer Linear Algebra Appl, 2005, 12: 809-826

\bibitem{pdv2009} Papy J M, De Lathauwer L, Van Huffel S. Exponential data fitting using multilinear algebra: the decimative case.
Journal of Chemometrics, 2009, 23: 341-351

\bibitem{q2015} Qi L. Hankel tensors: Associated Hankel matrices and Vandermonde decomposition. Commun Math Sci, 2015, 13: 113-125


\bibitem{RG} Ramana M, Goldman A J. Some geometric results in semidefinite programming. Journal of Global Optimization, 1995, 7: 33-50

\bibitem{x2015}  Xu C. Hankel tensors, Vandermonde tensors and their positivities, to appear in: Linear Algebra and Its Applications

\end{thebibliography}
\end{document}